\documentclass[10pt]{article}

\usepackage [applemac] {inputenc}
\usepackage[english,french]{babel}  
\usepackage[T1]{fontenc}

\usepackage{amssymb,amsmath,amsthm,dsfont}

\usepackage{xcolor}
\usepackage{geometry}      
\title{Note sur la continuité de la projection dans les facteurs de Host-Kra.}

\author{{\Large Jean-Fran\c cois Bertazzon\footnote{\texttt{jeffbertazzon@gmail.com}}}\\
 {\it \small Laboratoire d’Analyse, Topologie et Probabilit\'es,} \\ 
  {\it \small Aix-Marseille Universit\'e,}
   \\ {\it \small Avenue de l'escadrille Normandie-Ni\' emen. 13397  Marseille, France} }

\begin{document}

\makeatletter
\renewcommand\section{\@startsection{section}{1}{\z@}{6pt\@plus0pt}{6pt}{\scshape\centering}}
\renewcommand\subsection{\@startsection {subsection}{2}{\z@}{6pt \@plus 0ex  \@minus 0ex}{-6pt \@plus 0pt}{\reset@font\itshape \ \ \ \ \ \ \ \ }}
\makeatother

\makeatletter
\renewenvironment{thebibliography}[1]
     {\section*{\refname}%
      \@mkboth{\MakeUppercase\refname}{\MakeUppercase\refname}%
      \list{\@biblabel{\@arabic\c@enumiv}}%
           {\settowidth\labelwidth{\@biblabel{#1}}%
            \leftmargin\labelwidth
            \advance\leftmargin\labelsep
            \@openbib@code
            \usecounter{enumiv}%
            \let\p@enumiv\@empty
            \itemsep=0pt
            \parsep=0pt
            \leftmargin=\parindent
            \itemindent=-\parindent
            \renewcommand\theenumiv{\@arabic\c@enumiv}}%
      \sloppy
      \clubpenalty4000
      \@clubpenalty \clubpenalty
      \widowpenalty4000%
      \sfcode`\.\@m}
     {\def\@noitemerr
       {\@latex@warning{Empty `thebibliography' environment}}%
      \endlist}
\makeatother

\renewcommand{\proofname}{\textup{Preuve}}

\newcommand{\Ao}{\mathcal{A}} 
\newcommand{\Co}{\mathcal{C}} 
\newcommand{\Do}{\mathcal{D}}
\newcommand{\Xo}{\mathcal{X}}
\newcommand{\Zo}{\mathcal{Z}}

\newcommand{\EE}{\mathbb{E}} 

\newcommand{\bx}{\boldsymbol{x}}
\newcommand{\by}{\boldsymbol{y}}
\newcommand{\be}{{\boldsymbol{\varepsilon}}}
\newcommand{\bn}{\boldsymbol{n}}
\newcommand{\bm}{\boldsymbol{m}}
\newcommand{\m}{\overline{m}}
\newcommand{\x}{\overline{x}}

\newcommand{\e}{\varepsilon}
\newcommand{\Moy}{\frac{1}{N} \sum \limits_{n=0}^{N-1}}
\newcommand{\MoyN}{\frac{1}{N^\ell} \sum \limits_{n_1,\ldots,n_\ell=0}^{N-1}}
\newcommand{\Cl}{{V_\ell^*}}

\newtheorem{theoreme}{Théorème}
\newtheorem{proposition}{Proposition}
\newtheorem{corollaire}{Corollaire}

\maketitle

\renewcommand{\abstractname}{\textup{\textsc{Résumé}}}
\begin{abstract}
B. Host et B. Kra introduisent dans \cite{MR2150389} des facteurs caractéristiques pour l'étude de «moyennes ergodiques cubiques».
Ces facteurs permettent en particulier de résoudre des problèmes de récurrence multiple introduits par H. Furstenberg dans \cite{MR0498471}. 
Nous allons montrer que la continuité de la projection du système dans ses facteurs caractéristiques caractérisent la convergence des moyennes cubiques.
\end{abstract}

\section{Introduction et résultat.}

Soient $\ell\geq 1$ un entier et $\mathfrak X = (X,\Xo,\mu,T)$ un système dynamique ergodique.  Nous notons $\Cl=\{0,1\}^\ell\setminus\{(0,\ldots,0)\}$ et pour tous les  $\ell$-uplets $\be=(\e_1,\ldots,\e_\ell)$ et $\bn=(n_1,\ldots,n_\ell)$, nous notons $\be\cdot\bn = \e_1n_1+\cdots+\e_\ell n_\ell$. Nous nous intéressons à la convergence des moyennes cubiques
\begin{equation} \label{moy}
\MoyN \prod \limits_{\be\in\Cl} f_\be \big{(} T^{\bn \cdot \be} x \big{)},
\end{equation}
pour toutes les fonctions  $(f_\be)_{\be\in\Cl}$  mesurables et bornées.

Pour tout entier $\ell\geq2$, B. Host et B. Kra associent  dans \cite{MR2150389} au système ergodique $\mathfrak X$  un facteur $\mathfrak Z_\ell = (Z_\ell,\Zo_\ell,\mu_\ell,T_\ell)$. Ce facteur est une limite projective (topologique et en mesure) de nilsystèmes d'ordre inférieur à $\ell$. Nous renvoyons à \cite{MR2544760,HKM} pour la définition des nilsystèmes. Pour toute fonction $f\in$ L$^1(\mu)$, nous notons $\EE(f|\mathfrak Z_{\ell})$ l'espérance conditionnelle de la fonction $f$ par rapport à la tribu image réciproque de $\Zo_{\ell}$ par l'application facteur de $\mathfrak X$ dans $\mathfrak Z_\ell$. Ils montrent le résultat suivant :

\begin{theoreme}[B. Host \& B. Kra \cite{MR2150389}] \label{HK}
Soient $(X,\Xo,\mu,T)$ un système dynamique ergodique et $\ell\geq2$ un entier. Le facteur $\mathfrak Z_{\ell-1}$ est caractéristique pour la convergence L$^2$ des moyennes cubiques, dans le sens où pour toutes les fonctions mesurables et bornées $(f_\be)_{\be\in \Cl}$, les moyennes 
\[
\MoyN \prod \limits_{\be\in\Cl} f_\be \big{(} T^{\bn \cdot \be} x \big{)} \mbox{ et } \MoyN \prod \limits_{\be\in\Cl} \EE\big{(} f_\be | \mathfrak Z_{\ell-1}\big{)} \big{(} T^{\bn \cdot \be} x \big{)}
\]
convergent pour la norme L$^2(\mu)$ lorsque $N$ tend vers l'infini, et ont la même limite.
\end{theoreme}

De nombreux travaux concernent l'étude de la convergence des moyennes \eqref{moy}. Citons deux résultats que nous réutiliserons par la suite :

\begin{theoreme}[B. Host \& B. Kra \cite{MR2544760}]  \label{HKunif}
Soient $(X,\Xo,\mu,T)$ une limite projective topologique de nilsystèmes minimaux d'ordre fini et $\ell\geq2$ un entier. Pour toutes les fonctions continues  $(f_\be)_{\be\in \Cl}$, la suite de fonctions 
\[
\MoyN \prod \limits_{\be\in\Cl} f_\be \circ T^{\bn \cdot \be} \mbox{ converge uniformément lorsque $N$ tend vers l'infini.}
\]
\end{theoreme}

\begin{theoreme}[I. Assani \cite{MR2753294}]  \label{Assani}
Soient $(X,\Xo,\mu,T)$ un système dynamique mesuré, $\ell\geq2$ un entier et $(f_\be)_{\be\in \Cl}$ des fonctions mesurables et bornées sur $X$. Alors, les moyennes cubiques \eqref{moy} convergent presque sûrement. Nous noterons $\Do_\ell(f_\be:\be\in\Cl)$ la fonction limite.
\end{theoreme}

Si toutes les fonctions $f_\be$ sont égales à la même fonction $f$, nous désignons par $\Do_\ell(f)$ la $\ell$-ième fonction duale de $f$.

Contrairement aux moyennes ergodiques «de Furstenberg» (de la forme $(1/N)\sum _0^{N-1} f\circ T^n \cdots f\circ T^{kn}$), I. Assani remarqua dans \cite{assani} que la convergence des moyennes cubiques se généralisait au cas de transformations qui ne commutent pas (voir par exemple \cite{CF} pour une généralisation du théorème \ref{Assani}).

Soient $X$ un espace compact et $\mu$ une mesure borélienne sur $X$.  Alors une fonction $f$ sur $X$ est dite \textit{essentiellement continue} s’il existe une fonction  continue $f '$ sur $X$ égale à $f$ pour $\mu$-presque tout point.

Soient $(X, T)$ et $(Y,S)$ des systèmes topologiques, et $\mu$ et $\nu$ des mesures invariantes respectivement sur $X$ et $Y$.
Une application facteur au sens mesurable $\pi:X \mapsto Y$ est dite \textit{essentiellement continue} s’il existe une application facteur au sens topologique $\pi'$ telle que $\pi$ et $\pi'$ coïncident $\mu$-presque partout.

On remarque que dans ces définitions, on suppose seulement que $f$ et $\pi$ sont  déﬁnies presque partout. Si le support topologique de la mesure $\mu$ est égal à $X$, 
alors $f$ (resp. $\pi$) détermine $f'$ (resp. $\pi'$) de façon unique.

Nous montrons le résultat suivant.

\begin{theoreme} \label{th}
Soient $\mathfrak X=(X,d,T)$ un système dynamique topologique et $\ell$ un entier. Soit $\mu$ une mesure ergodique invariante de support topologique dense dans $X$.

Alors  la projection dans le $(\ell-1)$-ième facteur de Host-Kra est essentiellement continue si et seulement si pour toutes les fonctions mesurables et bornées $(f_\be)_{\be\in\Cl}$, les moyennes cubiques \eqref{moy} convergent presque partout vers une fonction limite essentiellement continue.
\end{theoreme}

Il y a de nombreuses raisons pour lesquelles il est intéressant de considérer des systèmes topologiques pour lesquels la projection dans ces facteurs sont continues.
Que ce soit pour des questions de convergence de moyennes ergodiques pondérées \cite{MR2544760}, pour l'étude de systèmes pour lesquels les applications duales  de toutes les fonctions continues sont continues \cite{HKM}, ou pour l'étude de suites universellement bonnes pour la convergence en moyenne de moyennes ergodiques multiples \cite{HKMbis}.

Notons que pour tout entier $\ell$, nous pouvons déduire directement du travail de B. Weiss \cite{MR799798,MR1727510} que tout système ergodique est conjugué en mesure à un système minimal et uniquement ergodique tel que la projection dans son $\ell$-ième facteur de Host-Kra est continue. 

Nous ne savons pas à l'heure actuelle si  les moyennes \eqref{moy} convergent partout (et uniformément) lorsque  le système est minimal et uniquement ergodique et que la projection dans le $(\ell-1)$-ième facteur de Host-Kra est continue.

\section{Preuve du théorème \ref{th}.}

\subsection{Supposons que la projection $\pi_{\ell-1}$ de $\mathfrak X$ dans $\mathfrak Z_{\ell-1}$ est continue.}

Fixons dans toute cette preuve $2^\ell-1$ fonctions mesurables et bornées $(f_\be)_{\be \in \Cl}$. Nous notons également $\tilde{f}_\be$ la fonction mesurable définie presque partout sur $Z_{\ell-1}$, telle que pour tout $\be\in \Cl$ et pour $\mu$-presque tout point $x$ de $X$, $\tilde{f}_\be\circ \pi_{\ell-1}(x) = \EE(f_\be|\mathfrak Z_{\ell-1})(x)$. Quitte à redéfinir ces fonctions sur des sous-ensembles dont la mesure est nulle, nous pouvons supposer qu'elles sont définies partout et qu'elles vérifient :
\[
\sup \left\{\left| \EE\big{(} f_\be|\mathfrak Z_{\ell-1}\big{)} (x)\right| ;  x \in X  \right\}= \sup \left\{ \left|\tilde{f}_\be (\zeta) \right|; \zeta \in Z_{\ell-1} \right \} =\sup \left\{\left|  f_\be (x)\right| ;  x \in X  \right\}.
\]

Comme l'a déjà observé I. Assani dans \cite{MR2753294}, on déduit directement des théorèmes \ref{HK} et \ref{Assani} que le facteur $\mathfrak Z_{\ell-1}$ est caractéristique pour la convergence \textit{presque sûre}  des moyennes cubiques de taille $\ell$, c'est-à-dire :
\begin{equation} \label{dual}
\Do_\ell(f_\be:\be\in \Cl) = \Do_\ell(\tilde{f}_\be;\be\in \Cl) \circ \pi_{\ell-1} \mbox{ pour $\mu$-presque tout point.}
\end{equation}

L'application $\pi_{\ell-1}$ est continue par hypothèse, il ne reste plus qu'à montrer que la fonction $\Do_\ell(\tilde{f}_\be;\be\in \Cl)$ est essentiellement continue. 

Pour tout entier $n$ et $\be\in \Cl$, on fixe une fonction continue $\tilde{f}^{(n)}_\be$ sur $Z_{\ell-1}$, bornée par $||\tilde{f}_\be||_\infty$ et égale à $\tilde{f}^{(n)}_\be$ à part peut-être sur un ensemble $\tilde{A}^{(n)}_\be$ de mesure inférieur à $1/2^{n+1}$.
Soit $Z_{\ell-1}^0$ l'ensemble des points génériques pour les fonctions $\mathds{1}(\tilde{A}^{(k)}_\be)$, et $n$ et $m$ deux entiers.

Pour tous les points de $Z_{\ell-1}^0$, puisque les fonctions $\tilde{f}^{(n)}_\be$  et $\tilde{f}^{(n+m)}_\be$ diffèrent sur un ensemble de taille au plus $2/2^{n+1}$, on trouve :
\[
\left|\Do_\ell \big{(} \tilde{f}^{(n+m)}_\be:\be\in\Cl\big{)}(x) - \Do_\ell \big{(} \tilde{f}^{(n)}_\be:\be\in\Cl\big{)}(x)\right| \leq \frac{1}{2^n} \prod \limits_{\be \in \Cl} || f_\be||_\infty.
\]

Mais puisque les fonctions $\Do_\ell  ( \tilde{f}^{(n)}_\be:\be\in\Cl )$ sont continues par le théorème \ref{HKunif}  et puisque la mesure $\mu$ a un support dense, la relation précédente se réécrit 
\[
\left |\left|\Do_\ell  ( \tilde{f}^{(n+m)}_\be:\be\in\Cl ) - \Do_\ell  ( \tilde{f}^{(n)}_\be:\be\in\Cl ) \right|\right|_\infty \leq \frac{1}{2^n} \prod \limits_{\be \in \Cl} || f_\be||_\infty.
\]

Il est alors clair que la fonction $\Do_\ell \big{(} \tilde{f}_\be:\be\in\Cl\big{)}$ est égale presque partout à la limite uniforme de la suite de fonctions $ \Do_\ell \big{(} \tilde{f}^{(n)}_\be:\be\in\Cl\big{)}$. C'est donc en particulier une fonction essentiellement continue.

\subsection{Supposons que pour toutes les fonctions mesurables et bornées $(f_\be)_{\be\in\Cl}$, les fonctions limites $\Do_\ell(f_\be:\be\in \Cl)$ sont essentiellement continues.}

Afin de simplifier les notations, nous supposons que l'application facteur $\pi_{\ell-1}$ de $\mathfrak X$ dans $\mathfrak Z_{\ell-1}$ est déﬁnie partout.

Choisissons une famille dénombrable de fonctions continues $(g_n)_n$ dense dans l'ensemble des fonctions continues $\Co(Z_{\ell-1} )$ pour la norme uniforme. Par hypothèse, chacune des fonctions $\Do_\ell (g_n \circ \pi_{\ell-1})$ est essentiellement continue sur $X$. Par l'équation \eqref{dual}, cette fonction est égale $\mu$-presque partout à $\Do_\ell (g_n) \circ \pi_{\ell-1}$ .

Il existe donc un sous-ensemble mesurable $X_0$ de mesure nulle pour $\mu$ tel que pour tout entier $n$, il existe donc une fonction continue $h_n$ sur $X$ telle que pour tout $x \in X \setminus X_0$, $\Do_\ell(g_n)\circ \pi_{\ell-1} (x) = h_n(x)$.

Mais la sous-algèbre engendrée par les fonctions $\Do_\ell(g)$, où $g$ est une fonction continue sur $Z_{\ell-1}$ est dense dans  $\Co(Z_{\ell-1})$ pour la norme uniforme \cite{MR2544760,HKM}. Donc pour toute fonction continue $g$ sur $Z_{\ell-1}$, il existe une fonction continue $h$ sur $X$ telle que 
\begin{equation} \label{un}
\mbox{pour tout point $x\in X\setminus X_0$, } g \circ \pi_{\ell-1}(x) = h(x).
\end{equation}

Puisque la mesure $\mu$ a un support dense dans $X$, l'ensemble $X_0$ est d'intérieur vide. En particuliers, deux fonctions continues qui coïncident sur l'ensemble $X_0$ sont égales. Plus précisemment si $f$ et $g$ sont deux fonctions continues sur $X$, $\sup\left\{ |f(x)-g(x)|;x\in X \setminus X_0\right\} = ||f-g||_\infty$. La fonction $h$ associée à $g$ en \eqref{un} est \textit{unique}, nous la notons $\Phi (g)$.

Par unicité, l'application $\Phi$ ainsi définie est compatible avec l'addition et la multiplication des fonctions. De plus, elle préserve la norme. On en déduit qu'il existe une application continue $\pi_{\ell-1}'$ de $X$ dans $Z_{\ell-1}$ telle que pour toute fonction continue $g$ sur $Z_{\ell-1}$, $\Phi( g) = g \circ \pi_{\ell-1}'$. On déduit de \eqref{dual} que $\pi_{\ell-1}=\pi_{\ell-1}'$ presque partout.

Il ne suffit alors plus qu'à vérifier que par unicité $\Phi(g\circ T) = \Phi(g) \circ T_{\ell-1}$ pour s'assurer que $\pi_{\ell-1}'$ est une application facteur topologique.
\qed

\section{Généralisation du théorème de Wiener-Wintner.}

Dans \cite{MR2544760}, B. Host et B. Kra généralisent le théorème de Wiener-Wintner de la manière suivante :

\begin{theoreme}[B. Host \& B. Kra \cite{MR2544760}] \label{thc}
Soient $\ell$ un entier et  $\mathfrak X=(X,d,T)$  un système  topologique minimal et uniquement ergodique tel que la projection dans son $\ell$-ième facteur de Host-Kra  est essentiellement continue.

Soient  $f_0$ une fonction continue sur $X$, $G$ un groupe de Lie nilpotent d'ordre inférieur à $\ell$, $\Gamma$ un sous-groupe co-compact de $G$, $g_0:G/\Gamma \to \mathbb{C}$ une fonction continue et  $y_0$ un point de $G$.

\begin{equation} \label{eq:HK}
\mbox{ Alors, la fonction } x \mapsto \lim \limits_{N \to \infty} \Moy f_0(T^n x)\cdot g_0(y_0^n \Gamma ) \mbox{ est définie partout}.
\end{equation}
\end{theoreme}

Nous proposons de préciser ce résultat de la manière suivante :

\begin{proposition} \label{pr} Sous les hypothèses du théorème \ref{thc}, la fonction définie en \eqref{eq:HK} est continue.
\end{proposition}

Nous ne savons pas si la continuité des fonctions limites \eqref{eq:HK} caractérise la continuité de la projection dans les facteurs de Host-Kra.

\begin{proof}[Preuve de la proposition \ref{pr}]
Nous nous plaçons sous les hypothèses du théorème \ref{thc}. Nous fixons une fonction continue $f_0$ définie sur $X$ ainsi qu'une fonction continue $g_0:G/\Gamma \to \mathbb{C}$, ou $G/\Gamma$ est une nilvariété d'ordre inférieur à $\ell$. Nous fixons de plus $y_0\in G$.

Nous définissons les fonctions $F_\infty$ et $F_N$ pour tout entier $N$, de $X$ dans $(\ell^{\infty} _\mathbb{C},||\cdot||_\infty)$ par :
\[
F_N(x) = \left( \Moy f_0(T^n  x) \cdot g_0(y_0^{n+m} \Gamma )  \right)_{m\in \mathbb{Z}} \mbox{ et }F_\infty(x) = \left( \lim_{N\to \infty } \Moy f_0(T^n  x)  \cdot g_0(y_0^{n+m} \Gamma )  \right)_{m\in \mathbb{Z}}.
\]

Si la fonction définie en \eqref{eq:HK} n'est pas continue en un certain point $x_\infty$, alors il existe $\e_0>0$  et une suite $(x_k)_k$ d'éléments de $X$ convergeant vers $x_\infty$ tels que pour tout entier $k$, $\left||  F_\infty( x_k) - F_\infty( x_\infty)  \right||_\infty \geq \e_0$.

Nous pouvons alors remarquer que pour  tout couple de points $(x,x')$ de $X^2$, la fonction $F_\infty$ vérifie
$||F_\infty(Tx)-F_\infty(Tx')||_\infty =||F_\infty(x)-F_\infty(x')||_\infty$. Mais puisque le point $x_\infty$ est transitif, alors la fonction limite $F_\infty$ n'est continue en aucun point.

D'autre part, puisque chaque fonction $F_N$ est continue, la fonction limite  $F_\infty$ admet par le théorème de la limite simple de Baire au moins un point de continuité ; ce qui est absurde.
\end{proof}

\textbf{Remerciements.} Je voudrais remercier B. Host pour son aide dans la préparation de cette note.

\end{document}